\documentclass[a4paper,10pt,leqno]{amsart}
\usepackage{amsmath,amsfonts,amssymb,amsthm,epsfig,epstopdf,url,array}
\usepackage{amsthm}
\usepackage{mathrsfs}
\usepackage{epsfig}
\usepackage{graphicx}

\newtheorem{thm}{Theorem}[section]
\newtheorem{prop}[thm]{Proposition}

\newtheorem{lem}[thm]{Lemma}
\newtheorem{rmk}[thm]{Remark}

\numberwithin{equation}{section}
\newcommand{\sgn}{\text{sgn}}
\newcommand{\ord}{\text{ord}}
\newcommand{\N}{{\mathbb{N}}}
\newcommand{\z}{{\mathbb{Z}}}
\newcommand{\q}{{\mathbb{Q}}}
\newcommand{\Tr}{\text{Tr}}
\newcommand{\qf}[1]{\langle #1 \rangle}
\newcommand{\conj}[1]{\overline{#1}}
\newcommand{\Case}[1]{\noindent\mbox{\bf Case #1}}
\newcommand{\comega}{\conj\omega}
\newcommand{\nequiv}{\not\equiv}
\newcommand{\gen}{\operatorname{gen}}
\newcommand{\disc}{\operatorname{disc}}
\newcommand{\rank}{\operatorname{rank}}
\renewcommand{\Re}{\operatorname{Re}}
\renewcommand{\Im}{\operatorname{Im}}
\newcommand{\mystrut}{\rule{0pt}{3ex}}
\newcommand{\ra}{\rightarrow}
\newcommand{\nra}{\not\rightarrow}
\newcommand{\floor}[1]{\left\lfloor #1 \right\rfloor}
\newcommand{\ceil}[1]{\left\lceil #1 \right\rceil}
\newcommand{\eq}{\equiv}
\newcommand{\rom}[1]{\uppercase\expandafter{\romannumeral #1\relax}}

\title[Almost regular $m$-gonal forms]{Representations of almost regular $m$-gonal forms $\rom2$}

\author{Dayoon Park}
\address{Department of Mathematical Sciences, Ulsan National Institute of Science and
Technology, Ulsan, Korea}
\email{pdy1016@unist.ac.kr}
\thanks {This work was supported by the National Research Foundation of Korea (NRF) grant funded by the Korea government(MSIT) (No. 2020R1A4A1016649). \\
This work was supported by the UBSI Fellowship Program (project No.1.210131.01) of UNIST}
\begin{document}

\maketitle

\begin{abstract}
It is known that any $m$-gonal form of $\rank n \ge 5$ is almost regular.
On the other words, any $m$-gonal form of $\rank n \ge 5$ represents every sufficiently large integer which is locally represented by the form.
    In this article, we treat the sufficiently large integers which are indeed represented by $m$-gonal forms of $\rank n \ge 5$ (those should be almost regular) among the locally represented integers.
\end{abstract}

\section{Introduction}
The {\it $m$-gonal number} is defined as a total number of dots to constitute a regular $m$-gon.
We especially call the total number of dots to constitute a regular $m$-gon with $x$ dots for each side
\begin{equation}\label{m number}
P_m(x)=\frac{m-2}{2}(x^2-x)+x
\end{equation}
as {\it $x$-th $m$-gonal number}.
Representing positive integers as a sum of polygonal numbers has been studied for a long time by many authors.
As a famous example, Fermat conjectured that every positive integer written as a sum of at most $m$ $m$-gonal numbers and which was resolved by Lagrange, Gauss, and Cauchy for $m=4$, $m=3$, and for all $m \ge 3$, respectively.

Following the definition of original $m$-gonal number, only positive integer $x$ would be admitted in (\ref{m number}).
But we may naturally generalize the $m$-gonal number by admitting every rational integer $x$ in (\ref{m number}).
The author completed the Fermat's conjecture of generalized $m$-gonal number virsion, i.e., classified the minimal $\ell_m$ for which every positive integer is written as a sum of $\ell_m$ generalized $m$-gonal numbers for all $m \ge 3$ in \cite{det}.
In order to more generally consider the problems about the representation as a sum of polygonal numbers, we define a weighted sum of $m$-gonal numbers 
\begin{equation} \label{m form}
F_m(\mathbf x)=a_1P_m(x_1)+\cdots+a_nP_m(x_n)
\end{equation}
as {\it $m$-gonal form} for $(a_1,\cdots,a_n) \in \N^n$.
For a positive integer $N \in \N$, if the diophantine equation $$F_m(\mathbf x)=N$$ has an integer solution $\mathbf x \in \z^n$, then
we say that {\it $F_m(\mathbf x)$ represents $N$}.
If an $m$-gonal form $F_m(\mathbf x)$ represents every positive integer, then we say that {\it $F_m(\mathbf x)$ is  universal}.
When the congruence equations
\begin{equation}\label{loc rep eq}F_m(\mathbf x) \eq N \pmod r\end{equation} have an integer solution $\mathbf x \in \z^n$ for all $r \in \z$, we say that {\it $F_m(\mathbf x)$ locally represents $N$}. 
In virtue of the Chinese Remainder Theorem, in discriminating the local representability, it would be enough to examine the solvability of the congruence equations \eqref {loc rep eq} modulo only the powers of primes.
In other words, we may claim that $F_m(\mathbf x)$ locally represents $N$ by showing that $F_m(\mathbf x)\equiv N \pmod{p^{\alpha}}$ has an integer solution $\mathbf x \in \z^n$ for every prime $p$ and every $\alpha \in \N$.
Which is equivalent with that $F_m(\mathbf x)=N$ has a $p$-adic integer solution $\mathbf x \in \z_p^n$ for every prime $p$.
Undoubtedly, if $F_m(\mathbf x)$ represents $N \in \N$, then $F_m(\mathbf x)$ locally represents $N \in \N$.
But the converse does not hold in general.
Especially, when the converse also holds (i.e.,  when an $m$-gonal form $F_m(\mathbf x)$ represents all positive integers which are locally represented), we say that the $m$-gonal form $F_m(\mathbf x)$ is {\it regular}.
The author and Kim completely determined every type of regular $m$-gonal form of $\rank \ge 4$ for all $m\ge 14$ with $m \not\eq 0 \pmod 4$ and $m\ge 28$ with $m \eq 0 \pmod 4$ in \cite{reg}.
If an $m$-gonal form $F_m(\mathbf x)$ represents all positive integers which are locally represented but finitely many, then we call $F_m(\mathbf x)$ is {\it almost regular}.
On the other words, almost regular $m$-gonal form $F_m(\mathbf x)$ would represent every sufficiently large integer which is locally represented by $F_m(\mathbf x)$ so we could take a sufficiently large constant $N_{m;(a_1,\cdots,a_n)}>0$ for which
$F_m(\mathbf x)=a_1P_m(x_1)+\cdots +a_nP_m(x_n)$ represents $N$ provided that
\begin{equation} \label{sl}
    \begin{cases}N \ge N_{m;(a_1,\cdots ,a_n)} \\ N \text{ is locally represented by $F_m(\mathbf x)$}. \end{cases}
\end{equation}    
Some authors have made an effort to examine a  constant $N_{m;(a_1,\cdots,a_n)}$ in \eqref{sl} for some specific type of almost regular $m$-gonal forms $a_1P_m(x_1)+\cdots +a_nP_m(x_n)$ (see \cite{N1}, \cite{N2}, \cite{MS}, and \cite{K}).

\vskip 1em

By Theorem 4.9 (1) in \cite{CO}, every quadratic polynomial of $\rank \ge 5$ whose quadratic part is positive definite is almost regular.
Which says that for any $m$-gonal form $F_m(\mathbf x)=a_1P_m(x_1)+\cdots+a_nP_m(x_n)$ of rank $n \ge 5$, we could always take a sufficiently large constant $N_{m;(a_1,\cdots,a_n)}$ for which $F_m(\mathbf x)$ represents $N$ satisfying \eqref{sl}.
Indeed, the $5$ is optimal, i.e., one may find an $m$-gonal form of $\rank 4$ which is not almost regular.
In this article, we more generally ponder the behavior pattern of sufficiently large integers $N_{m;(a_1,\cdots,a_n)}$ 
for arbitrary $(a_1,\cdots,a_n) \in \N^n$ with $n \ge 5$.
In \cite{rank 6}, the author firstly looked into
$N_{m;(a_1,\cdots,a_n)}$ when $n \ge 6$ leaving the case $n=5$.
In this article, we completely resolve the issue for all rank $n \ge 5$.
The following theorem is the main theorem in this article.

\begin{thm}\label{main thm}
Let $F_m(\mathbf x)=a_1P_m(x_1)+\cdots+a_nP_m(x_n)$ be an $m$-gonal form of rank $n \ge 5$.
Then $F_m(\mathbf x)$ represents any integer $N$ provided that
\begin{equation} \label{main eq}
    \begin{cases}N \ge N(a_1,\cdots,a_n)\cdot (m-2)^3 \\ N \text{ is locally represented by } F_m(\mathbf x)\end{cases}
\end{equation}
where $N(a_1,\cdots,a_n)$ is a constant which is dependent only on $a_1,\cdots,a_n$. 
\end{thm}

In this sense, both of the cubic with respect to $m$ in (\ref{main eq}) and the rank $n\ge5$ are optimal.
In order to prove Theorem \ref{main thm}, we use the arithmetic theory of quadratic form.

\vskip 0.3em

\begin{rmk}
In \cite{KL}, Kane and Liu proved that the asymptotically increasing $\gamma_m$ ($\gamma_m$ is the smallest positive integer for which an $m$-gonal form is universal if the $m$-gonal form represents every positive integer up to $\gamma_m$) is bounded by $C_{\epsilon}\cdot m^{7+\epsilon}$ for any $\epsilon >0$ and an absolute constant $C_{\epsilon}>0$ which is dependent on $\epsilon$ by observing the sufficiently large integers which are represented by the $m$-gonal forms of escalator tree of $m$-gonal forms.
Theorem \ref{main thm} may directly improve the upper bound for $\gamma_m$ as $C\cdot (m-2)^3$ for some absolute constant $C$ since every $m$-gonal form of escalator tree of $m$-gonal forms of depth $5$ is locally universal.
Although, the exact growth of $\gamma_m$ on $m$ was already revealed as linear on $m$ by the author and Kim \cite{KP'}.
\end{rmk}

\section{preliminaries}

Before we move on, we set our languages and notations which are used throughout this article.
Basically, our languages and notations follow the convention in the arithmetic theory of quadratic form.
A {\it quadratic form over $\z$ (or $\z_p$)} is a $\z$ (or $\z_p$)-module $L$ equipped with a quadratic form $Q:L\rightarrow \z$ (or $\z_p$) and the  corresponding symmetric bilinear map $B(x,y):=\frac{1}{2}(Q(x+y)-Q(x)-Q(y))$.
For a quadratic form $L=\z v_1+\cdots+\z v_n$ of $\rank n$, we usually identify $L$ with its {\it Gram-matrix $\begin{pmatrix}B(v_i,v_j)\end{pmatrix}_{n\times n}$} with respect to a basis $\{v_1,\cdots,v_n\}$ of $L$.
When $L$ is a binary (i.e., $n=2$) quadratic form over $\z$, we adopt the notation $[Q(v_1),B(v_1,v_2),Q(v_2)]$ instead of $\begin{pmatrix}B(v_i,v_j)\end{pmatrix}_{2 \times 2}$.
We write an $m$-gonal form 
\begin{equation} \label{m form}
a_1P_m(x_1)+\cdots+a_nP_m(x_n)
\end{equation}
 simply as $\left<a_1,\cdots,a_n\right>_m$.
When $m=4$, we conventionally denote the square form (i.e., diagonal quadratic form) of (\ref{m form}) by $\left<a_1,\cdots, a_n\right>$ instead of $\left<a_1,\cdots,a_n\right>_4$.

\vskip 0.3em

Note that an integer $A(m-2)+B$ is represented by an $m$-gonal form 
$$\sum_{i=1}^na_iP_m(x_i)=\frac{m-2}{2}\{(a_1x_1^2+\cdots+a_nx_n^2)-(a_1x_1+\cdots+a_nx_n)\}+(a_1x_1+\cdots+a_nx_n)$$
over $\z$ (or $\z_p$)
if and only if the diophantine system
\begin{equation}\label{eq1}
\begin{cases}a_1x_1^2+\cdots+a_nx_n^2=2A+B+k(m-4) \\ a_1x_1+\cdots+a_nx_n=B+k(m-2)  \end{cases}\end{equation}
is solvable for some $k \in \z$ (or $\z_p$) over $\z$ (or $\z_p$).
Emphatically note that the system (\ref{eq1}) holds if and only if the equation
\begin{equation} \label{eq2}
    \left(B+k(m-2)-\left(\sum \limits_{i=2}^na_ix_i\right)\right)^2+\sum\limits_{i=2}^na_1a_ix_i^2=2Aa_1+Ba_1+k(m-4)a_1 
\end{equation}
holds with $x_1=\frac{1}{a_1}\left(B+k(m-2)-\left(\sum \limits_{i=2}^na_ix_i\right)\right)$.
The left hand side of the equation (\ref{eq2}) may be written as
$$Q_{a_1 ; \mathbf a}((x_2,\cdots,x_n)-(B+k(m-2))(r_2,\cdots,r_n)) + (B+k(m-2))^2 \cdot \left(1-\sum \limits _{i=2}^na_ir_i\right) $$
where $Q_{a_1 ; \mathbf a}(x_2,\cdots,x_n):=\sum_{i=2}^n(a_1a_i+a_i^2)x_i^2+\sum_{2\le i<j \le n}2a_ia_jx_ix_j$ is a quadratic form and $r_2,\cdots,r_n \in \q$ are the solution for
$$\begin{cases}
(a_1a_2+a_2^2)r_2+a_2a_3r_3+\cdots+a_2a_nr_n=a_2\\
a_2a_3r_2+(a_1a_3+a_3^2)r_3 +\cdots + a_3a_nr_n=a_3 \\
\quad \quad \quad \quad \quad \quad \quad \quad \quad \vdots \\
a_2a_nr_2+a_3a_nr_3+\cdots +(a_1a_n+a_n^2)r_n=a_n
\end{cases}$$
(in practice, $r_2=\cdots=r_n=\frac{1}{a_1+\cdots+a_n}$). 
Note that since 
\begin{equation}\label{n-1 sub}
Q_{a_1 ; \mathbf a}= \begin{pmatrix}
a_1a_2+a_2^2 & a_2a_3 &\cdots & a_2a_n\\
a_2a_3 & a_1a_3+a_3^2 & \cdots & a_3a_n\\
\vdots & \vdots & \ddots  & \vdots \\
a_2a_n & a_3a_n & \cdots &a_1a_n+a_n^2
\end{pmatrix}
\end{equation}
is a sub-quadratic form of rank $n-1$ of the positive definite diagonal quadratic form 
$\sum_{i=1}^na_1a_ix_i^2$ with $a_1x_1+\cdots+a_nx_n=0$,
we have that the determinant of $Q_{a_1 ; \mathbf a}$ is positive.
We denote the non-zero determinant of $Q_{a_1 ; \mathbf a}$ as $d_{a_1;\mathbf a}$.

\vskip 1em

Our main consideration in this article is for an integer $A(m-2)+B$ which is locally represented by an $m$-gonal form $\left<a_1,\cdots,a_n\right>_m$ of rank $n \ge 5$, the $k \in \z$ for which the equation \eqref{eq2} has a local solution of bounded order over local ring $\z_p$ for every prime $p$.

\vskip 0.3em

In this article, to avoid the appearance of superfluous duplications, we always assume that an $m$-gonal form $\left<a_1,\cdots,a_n\right>_m$ is primitive, i.e., $\gcd(a_1,\cdots,a_n)=1$.
Any unexplained terminology and notation can be found in \cite{O}.

\section{Local representation}

The following proposition makes determine the local representability of an integer by an $m$-gonal form.

\begin{prop} \label{loc.rep}
Let $F_m(\mathbf x)=a_1P_m(x_1)+\cdots+a_nP_m(x_n)$ be a primitive $m$-gonal form.
\begin{itemize}
    \item [(1) ] When $p$ is an odd prime with $p|m-2$, $F_m(\mathbf x)$ is universal over $\z_p$.
    \item [(2) ] When $m \not\eq 0 \pmod 4$, $F_m(\mathbf x)$ is universal over $\z_2$.
   \item [(3) ] When $p$ is an odd prime with $(p,m-2)=1$, an integer $N$ is represented by $F_m(\mathbf x)$ over $\z_p$ if and only if the integer $8(m-2)N+(a_1+\cdots+a_n)(m-4)^2$ is represented by the diagonal quadratic form $\left<a_1,\cdots,a_n \right>$ over $\z_p$.
    \item [(4) ] When $m \eq 0 \pmod 4$, an integer $N$ is represented by $F_m(\mathbf x)$ over $\z_2$ if and only if the integer $\frac{m-2}{2}N+(a_1+\cdots+a_n)\left(\frac{m-4}{4}\right)^2$ is represented by the diagonal quadratic form $\left<a_1,\cdots,a_n \right>$ over $\z_2$.
\end{itemize}
\end{prop}
\begin{proof}
(1) One may easily show that for any $N \in \z_p$, $P_m(x)=N$ is solovable over $\z_p$ by using Hensel's Lemma (13:9 in \cite{O}).
Which yields the claim from the assumption that $F_m(\mathbf x)$ is primitive.\\
(2) For any $N \in \z_2$, the equation $P_m(x)=N$ has a $2$-adic integer solution $\frac{m-4-\sqrt{(m-4)^2-8N(m-2)}}{2(m-2)} \in \z_2$.
Which yields the claim from the assumption that $F_m(\mathbf x)$ is primitive. \\
(3) Note that $$N=\sum_{i=1}^na_iP_m(x_i)$$
if and only if
$$8(m-2)N+(a_1+\cdots+a_n)(m-4)^2=\sum_{i=1}^na_i(2(m-2)x_i-(m-4))^2$$
where $2(m-2) \in \z_p^{\times}$.
It completes the proof.\\
(4) Note that $$N=\sum_{i=1}^na_iP_m(x_i)$$
if and only if
$$\frac{m-2}{2}N+(a_1+\cdots+a_n)\left(\frac{m-4}{4}\right)^2=\sum_{i=1}^na_i\left(\frac{m-2}{2}x_i-\frac{m-4}{4}\right)^2$$
where $\frac{m-2}{2} \in \z_2^{\times}$ and $\frac{m-4}{4} \in \z_2$.
It completes the proof.
\end{proof}

\vskip 0.8em

\begin{prop}\label{odd,5}
For an odd prime $p$, if there are more than or equal to five units of $\z_p$ in $\{a_1,\cdots,a_n\}$  by admitting a recursion, then
for any $A,B,k \in \z_p$,
the equation \eqref{eq2} has a primitive solution $(x_2,\cdots,x_n) \in \z_p^{n-1}$.
\end{prop}
\begin{proof}
When $a_1+\cdots+a_n \in \z_p^{\times}$, the quadratic form $Q_{a_1;\mathbf a}\otimes \z_p(=\left<a_1+\cdots+a_n\right>^{\perp})$ contains a quaternary unimodular $\z_p$-lattice which is universal and isotropic and $r_2=\cdots=r_n=\frac{1}{a_1+\cdots+a_n}\in \z_p$.
So we may take a $(x_2,\cdots,x_n) \in \z_p^{n-1}$ satisfying \eqref{eq2} with primitive $\mathbf x-(B+k(m-2))\mathbf r \in \z_p^{n-1}$.
If $(B+k(m-2))\mathbf r$ is non-primitive, we obtain our claim.
If $(B+k(m-2))\mathbf r$ is primitive, by changing $\mathbf x$ to $-\mathbf x+2(B+k(m-2))\mathbf r$ if it is necessary, we obtain our claim. 

Now assume that $a_1+\cdots+a_n \in p\z_p$.
Recall that the equation \eqref{eq2} holds with $a_1x_1=a_2x_x+\cdots+a_nx_n-(B+k(m-2))$
if and only if the system \eqref{eq1} holds.
Without loss of generality, let $a_2 \in \z_p^{\times}$.
And then one may see that the equation \eqref{eq2} has a primitive solution $(x_2,\cdots,x_n) \in \z_p^{n-1}$ with $x_2=0$ similarly with the above since a ternary unimodular $\z_p$-lattice is also universal and isotropic.
\end{proof}

\vskip 0.8em

\begin{rmk}
Recall that $A(m-2)+B$ is represented by 
$\left<a_1,\cdots,a_n\right>_m$ over $\z_p$ if and only if there is $k \in \z_p$ for which
$$\begin{cases}a_1x_1^2+\cdots+a_nx_n^2=2A+B+k(m-4) \\ a_1x_1+\cdots+a_nx_n=B+k(m-2)  \end{cases}$$
is solvable over $\z_p$ and which is equivalent with that the equation 
$$\left(B+k(m-2)-\left(\sum \limits_{i=2}^na_ix_i\right)\right)^2+\sum\limits_{i=2}^na_1a_ix_i^2=2Aa_1+Ba_1+k(m-4)a_1 $$ is solvable over $\z_p$.
By Proposition \ref{odd,5}, when an integer $A(m-2)+B$ is locally represented by  $\left<a_1,\cdots,a_n\right>_m$ of rank $n \ge 5$, for any integer $k \in \z$, the above equation is solvable (even primitively solvable) over $\z_p$ 
for almost all primes $p$.
Now, in the rest of this section, we consider the $k \in \z_p$ over $\z_p$ for finitely many remaining primes $p$.
\end{rmk}

\vskip 0.8em

\begin{prop} \label{primitive}
Let an integer $N$ is represented by $m$-gonal form $\left<a_1,\cdots,a_n\right>_m$ of rank $n \ge 5$ over a non-dyadic $\z_p$.
Then there is a vector $\mathbf x \in \z_p^n$ for $$a_1P_m(x_1)+\cdots+a_nP_m(x_n)=N$$
such that $(x_2,\cdots,x_n) \in \z_p^{n-1}$ is primitive otherwise $p|m-4$ and $\ord_p(a_1) \ge \max\{\ord_p(a_2), \cdots, \ord_p(a_n)\}$ with anisotropic quadratic form $\left<a_2,\cdots,a_n\right>\otimes \z_p$.
\end{prop}
\begin{proof}
First, assume that $p|m-2$.
Without loss of generality, we may assume that $\left<a_1,\cdots,a_{n-1}\right>_m$ is primitive (by changing the order of $a_2,\cdots,a_n$ if it is necessary).
By Proposition \ref{loc.rep}, since $\left<a_1,\cdots,a_{n-1}\right>_m$ is universal over $\z_p$, there would be $(x_1,\cdots,x_{n-1}) \in \z_p^{n-1}$ for which
$$a_1P_m(x_1)+\cdots+a_{n-1}P_m(x_{n-1})+a_nP_m(1)=N$$
with primitive $(x_1,\cdots,x_{n-1},1) \in \z_p^{n-1}$.

Secondly, assume that $(p,m-2)=(p,m-4)=1$. 
Note that $$P_m(x)=P_m\left(-x+(m-4)(m-2)^{-1}\right).$$
For any $\mathbf x \in \z_p^n$ with $a_1P_m(x_1)+\cdots+a_nP_m(x_n)=N$, by changing $\mathbf x$ to $-\mathbf x+(m-4)(m-2)^{-1}$ if it is necessary, we obtain the claim.

Finally, assume that $p|m-4$. 
Note that 
$$a_1P_m(x_1)+\cdots +a_nP_m(x_n)=N$$
if and only if
$$\sum_{i=1}^n a_i(2(m-2)x_i-(m-4))^2=8(m-2)N+(m-4)^2(a_1+\cdots+a_n).$$
If the quadratic form $\left<a_2,\cdots,a_n\right>\otimes \z_p$ is isotropic, then we are done.
If not, since $\left<a_1,\cdots,a_n\right>\otimes \z_p$ with $n\ge 5$ is isotropic, there would be a primitive $\mathbf x \in \z_p^n$ for $\sum_{i=1}^n a_i(2(m-2)x_i-(m-4))^2=8(m-2)N+(m-4)^2(a_1+\cdots+a_n)$.
If $(x_2,\cdots,x_n)$ is primitive, then we are done.
Otherwise, we have that $x_1 \in \z_p^{\times}$.
From our assumption, we may assume that  $\ord_p(a_1) < \ord_p(a_2)$ by changing the order of $a_2,\cdots,a_n$ if it is necessary.
We could take $x_1' \in \z_p$ such that 
$a_1(2(m-2)x_1-(m-4))^2+a_2(2(m-2)x_2-(m-4))^2=a_1(2(m-2)x_1'-(m-4))^2+a_2(2(m-2)\cdot 1-(m-4))^2$ since $(u+p\z_p)^2=u^2+p\z_p$ for any $u \in \z_p^{\times}$.
And then we may obtain that $a_1P_m(x_1')+a_2P_m(1)+a_3P_m(x_3)+\cdots+a_nP_m(x_n)=N$ holds with primitive $(1,x_3,\cdots,x_n) \in \z_p^{n-1}$.
\end{proof}

\vskip 0.8em

\begin{rmk}
By Proposition \ref{primitive}, for an integer $A(m-2)+B$ which is locally represented by an $m$-gonal form $\left<a_1,\cdots,a_n\right>_m$ of rank $n \ge 5$, there is $k \in \z_p$ for which the equation \eqref{eq2} is not only solvable but also primitively solvable over $\z_p$
when $p$ is an odd prime otherwise $p|m-4$ and $\ord_p(a_1) \ge \max\{\ord_p(a_2), \cdots, \ord_p(a_n)\}$ with anisotropic quadratic form $\left<a_2,\cdots,a_n\right>\otimes \z_p$.
In the following proposition, we describe the distribution of such $k \in \z_p$.
\end{rmk}

\vskip 0.8em

\begin{prop} \label{odd,nd12}
If the equation \eqref{eq2} is primitively solvable over non-dyadic $\z_p$ for some $k \in \z_p$, then for any $p$-adic integer $k' \eq k \pmod{p^{1+2\ord_p(d_{a_1;\mathbf a})}}$, 
$$ \left(B+k'(m-2)-\left(\sum \limits_{i=2}^na_ix_i\right)\right)^2+\sum\limits_{i=2}^na_1a_ix_i^2=2Aa_1+Ba_1+k'(m-4)a_1$$
is also primitively solvable over $\z_p$.
\end{prop}
\begin{proof}
In virtue of the {\it Jordan decomposition} (see $\S$ 91 in \cite{O}), we may take $T \in GL_{n-1}(\z_p)$ for which
$$T^t\cdot  \begin{pmatrix}
a_1a_2+a_2^2 & a_2a_3 &\cdots & a_2a_n\\
a_2a_3 & a_1a_3+a_3^2 & \cdots & a_3a_n\\
\vdots & \vdots & \ddots  & \vdots \\
a_2a_n & a_3a_n & \cdots &a_1a_n+a_n^2
\end{pmatrix} \cdot T=\begin{pmatrix}
a_2' & 0 &\cdots & 0\\
0 & a_3' &\cdots & 0\\
\vdots & \vdots & \ddots  & \vdots \\
0 & 0 & \cdots &a_n'
\end{pmatrix}.$$
Which allows us to consider the simple diagonal form
$$\sum_{i=2}^na_i'(x_i-(B+k(m-2))r_i')^2$$
as a substitute for $Q_{a_1 ; \mathbf a}(\mathbf x-(B+k(m-2))\mathbf r)$ where $\mathbf r'=T\cdot \mathbf r$.
On the other words, we show that if
$$\sum_{i=2}^na_i'(x_i-(B+k(m-2))r_i')^2=(2A+B+k(m-4))a_1-(B+k(m-2))^2\left(1-\sum_{i=2}^na_ir_i\right)$$
is primitively solvable over $\z_p$, then
$$\sum_{i=2}^na_i'(x_i-(B+k'(m-2))r_i')^2=(2A+B+k'(m-4))a_1-(B+k'(m-2))^2\left(1-\sum_{i=2}^na_ir_i\right)$$
is also primitively solvable over $\z_p$ for any $k' \eq k \pmod{p^{1+2\ord_p(d_{a_1;\mathbf a})}}$.

Note that
\begin{equation}\label{ord r'a'}
\begin{cases}-\ord_p(d_{a_1;\mathbf a}) \le \ord_p(r_i') \\ 
0 \le \ord_p(a_i') \le \ord_p(d_{a_1;\mathbf a}) \end{cases}    
\end{equation}
for all $2 \le i \le n$,
\begin{equation} \label{square}
  \begin{cases}
(u+p\z_p)^2=u^2+p\z_p \\
\left(\z_p^{\times}+\frac{u}{p^r}\right)^2=\left(\frac{p^r\z_p^{\times}+u}{p^r}\right)^2=\frac{p^r\z_p^{\times}+u^2}{p^{2r}} \\
\left(\z_p+\frac{u}{p^r}\right)^2=\left(\frac{p^r\z_p+u}{p^r}\right)^2=\frac{p^r\z_p+u^2}{p^{2r}}
\end{cases}  
\end{equation}
for $u \in \z_p^{\times}$ and $r \in \z_{>0}$ and, $(2A+B+k(m-4))a_1-(B+k(m-2))^2\left(1-\sum_{i=2}^na_ir_i\right)\eq(2A+B+k'(m-4))a_1-(B+k'(m-2))^2\left(1-\sum_{i=2}^na_ir_i\right) \pmod{p^{1+\ord_p(d_{a_1;\mathbf a})}\z_p}$.
So if there is an $i$ such that $\ord_p((B+k(m-2))r_i')<0$, then we may obtain the claim from the second and third observations of \eqref{square}.

Now assume that $\ord_p((B+k(m-2))r_i') \ge 0$ for all $2\le i \le n$.
If $\ord_p((B+k(m-2))r_i')>0$ for all $2 \le i \le n$, then from the first observation of \eqref{square}, we obtain the claim.
Now we assume that $\ord_p((B+k(m-2))r_i')=0$ for some $2 \le i \le n$.
When $$\ord_p\left((2A+B+k(m-4))a_1-(B+k(m-2))^2\left(1-\sum_{i=2}^na_ir_i\right)\right)\le \ord_p(d_{a_1;\mathbf a}),$$
$(2A+B+k'(m-4))a_1-(B+k'(m-2))^2\left(1-\sum_{i=2}^na_ir_i\right)$ is equivalent with $(2A+B+k(m-4))a_1-(B+k(m-2))^2\left(1-\sum_{i=2}^na_ir_i\right)$ up to $\z_p$-unit square and the diagonal quadratic form $\left<a_2,\cdots,a_n\right>\otimes \z_p$ of rank $n-1 \ge 4$ is $p^{\max \limits _{2\le i \le n}\{\ord_p(a_i')\}}\z_p$-universal.
Therefore $(2A+B+k'(m-4))a_1-(B+k'(m-2))^2\left(1-\sum_{i=2}^na_ir_i\right)$ would be clearly represented by the diagonal quadratic form $\left<a_2,\cdots,a_n\right>\otimes \z_p$.
The primitively representability is obtained
by changing $x_i$ to $-x_i+2(B+k'(m-2))r_i'$ if it is necessary.
\end{proof}

\vskip 0.8em

\begin{rmk}
By Propositions \ref{primitive}, and \ref{odd,nd12}, we may obtain that when an integer $A(m-2)+B$ is locally represented by an $m$-gonal form $\left<a_1,\cdots,a_n\right>_m$
of rank $n \ge 5$, there is a residue $k_{p;\mathbf a}(A(m-2)+B)$ in $\z/p^{1+2\ord_p(d_{a_1;\mathbf a})}\z$ for which 
$$\text{the equation \eqref{eq2} is primitively solvable over non-dyadic local ring $\z_p$}$$
for any $k \eq k_{p;\mathbf a}(A(m-2)+B) \pmod{p^{1+2\ord_p(d_{a_1;\mathbf a})}}$ otherwise
\begin{equation} \label{bad}
\begin{cases}
p|m-4\\
\text{$\ord_p(a_1) \ge \max \limits_{2\le i \le n}\{\ord_p(a_i)\}$ with anisotropic $\left<a_2,\cdots,a_n\right>\otimes \z_p$.}
\end{cases}    
\end{equation}

On the other hand, there would be a chance that for an integer $N$ which is represented by $\left<a_1,\cdots,a_n\right>_m$ over $\z_p$,
for any $\mathbf x \in \z_p^n$ with $N=a_1P_m(x_1)+\cdots+a_nP_m(x_n)$, $(x_2,\cdots,x_n) \in \z_p^{n-1}$ is non-primitive (i.e., $\min\{\ord_p(x_2),\cdots,\ord_p(x_n)\}$ $>0$) when $p$ is an odd prime satisfying \eqref{bad}.

Nevertheless, it is a lucky misfortune that at least it would be possible to bound $\min \limits _{2\le i \le n}\{\ord_p(x_i)\}$
as smaller than or equal to $\frac{1}{2}\ord_p(a_1)$.
We see the claim in Proposition \ref{min ord}.
\end{rmk}

\vskip 0.8em

\begin{prop}\label{min ord}
Let $\left<a_1,\cdots,a_n\right>_m$ be an $m$-gonal form of rank $n \ge 5$ and $p$ be an odd prime satisfying \eqref{bad}.
For any integer $N$ which is represented by $\left<a_1,\cdots,a_n\right>_m$ over $\z_p$,
there is an $\mathbf x \in \z_p^n$ such that
$$\begin{cases}
a_1P_m(x_1)+\cdots+a_nP_m(x_n)=N \\
\min \{\ord_p(x_2), \cdots, \ord_p(x_n)\}\le \frac{1}{2}\ord_p(a_1).
\end{cases}$$
Besides, when the equation \eqref{eq2} has a solution with
$\min \limits _{2\le i \le n}\{\ord_p(x_i)\}\le \frac{1}{2}\ord_p(a_1)$,
for any $k' \eq k \pmod{p^{1+\ord_p(a_1)+2\ord_p(d_{a_1;\mathbf a})}}$, the equation
$$ \left(B+k'(m-2)-\left(\sum \limits_{i=2}^na_ix_i\right)\right)^2+\sum\limits_{i=2}^na_1a_ix_i^2=2Aa_1+Ba_1+k'(m-4)a_1$$
also has a solution with $\min \limits _{2\le i \le n}\{\ord_p(x_i)\}\le \frac{1}{2}\ord_p(a_1)$.
\end{prop}
\begin{proof}
Recall that
$$N=a_1P_m(x_1)+\cdots+a_nP_m(x_n)$$
if and only if
$$8(m-2)N+(m-4)^2(a_1+\cdots+a_n)=\sum_{i=1}^na_i(2(m-2)x_i-(m-4))^2.$$
When $\ord_p(8(m-2)N+(m-4)^2(a_1+\cdots+a_n)) < \ord_p(a_1)$, we directly obtain the claim for $\ord_p(a_i(2(m-2)x_i-(m-4))^2)$ should be less than or equal to $\ord_p(8(m-2)N+(m-4)^2(a_1+\cdots+a_n))$ for some $2\le i \le n$ (by changing $x_i$ to $-x_i+(m-4)(m-2)^{-1}$ if it is necessary).
Now assume that $\ord_p(8(m-2)N+(m-4)^2(a_1+\cdots+a_n)) \ge  \ord_p(a_1)$.
We may take $x_1 \in \z_p$ for
$$\ord_p(8(m-2)N+(m-4)^2(a_1+\cdots+a_n)-a_1(2(m-2)x_1-(m-4))^2)=\ord_p(a_1).$$
Since $\left<a_2,\cdots,a_n\right>\otimes \z_p$ of rank $n-1\ge 4$ is $p^{\max\limits_{2 \le i \le n}\{\ord_p(a_i)\}}\z_p$-universal and $\ord_p(a_1) \ge \max\limits_{2 \le i \le n}\{\ord_p(a_i)\}$, we may take $(x_2,\cdots,x_n) \in \z_p^{n-1}$
for
$$8(m-2)N+(m-4)^2(a_1+\cdots+a_n)=\sum_{i=1}^na_i(2(m-2)x_i-(m-4))^2.$$
And in this case too, we have  $min\{\ord_p(x_2),\cdots,\ord_p(x_n)\} \le \frac{1}{2}\ord_p(a_1)$ by the same token with the above.

For the second statement, one may show it through similar arguments with the proof of Proposition \ref{odd,nd12}.
Note that $\frac{1}{a_1+\cdots+a_n}=r_2=\cdots=r_n \in \z_p$.
\end{proof}

\vskip 0.8em

At the end of this section, we see similar argument with the above over the dyadic local ring $\z_2$ in Proposition \ref{2,nd}.
 
\begin{prop}\label{2,nd}
Let $\left<a_1,\cdots,a_n\right>_m$ be an $m$-gonal form of rank $n \ge 5$.
For any integer $N$ which is represented by $\left<a_1,\cdots,a_n\right>_m$ over $\z_2$,
there is an $\mathbf x \in \z_2^n$ such that
$$\begin{cases}
a_1P_m(x_1)+\cdots+a_nP_m(x_n)=N \\
\min \{\ord_2(x_2), \cdots, \ord_2(x_n)\}\le \frac{1}{2}\ord_2(a_1)+1.
\end{cases}$$
Besides, when the equation \eqref{eq2} has a solution with
$\min \limits _{2\le i \le n}\{\ord_2(x_i)\}\le \frac{1}{2}\ord_2(a_1)+1$ over $\z_2$,
for any $k' \eq k \pmod{2^{3+\ord_2(a_1)+2\ord_2(d_{a_1;\mathbf a})}}$,
$$ \left(B+k'(m-2)-\left(\sum \limits_{i=2}^na_ix_i\right)\right)^2+\sum\limits_{i=2}^na_1a_ix_i^2=2Aa_1+Ba_1+k'(m-4)a_1$$
also has a solution with $\min \limits _{2\le i \le n}\{\ord_2(x_i)\}\le \frac{1}{2}\ord_2(a_1)+1$ over $\z_2$.
\end{prop}
\begin{proof}
For the first statement, one may show this through similar processings with the proof of Propositions \ref{primitive} and \ref{min ord}.

For the second statement, one may show this through similar processing with the proof of Proposition \ref{odd,nd12} after Jordan decomposing (see $\S$ 91 in \cite{O}) $$\begin{pmatrix} 
a_1a_2+a_2^2 & a_2a_3 &\cdots & a_2a_n\\
a_2a_3 & a_1a_3+a_3^2 & \cdots & a_3a_n\\
\vdots & \vdots & \ddots  & \vdots \\
a_2a_n & a_3a_n & \cdots &a_1a_n+a_n^2
\end{pmatrix}\otimes \z_2.$$
One may use that
$$\begin{cases}
\begin{pmatrix}
0 & 1 \\ 1 & 0
\end{pmatrix}
=

\begin{pmatrix}
\frac{1}{2} & \frac{1}{2} \\ \frac{1}{2} &- \frac{1}{2}
\end{pmatrix}^t

\begin{pmatrix}
2 & 0 \\ 0 & -2
\end{pmatrix}

\begin{pmatrix}
\frac{1}{2} & \frac{1}{2} \\ \frac{1}{2} &- \frac{1}{2}
\end{pmatrix}

\\

\begin{pmatrix}
2 & 1 \\ 1 & 2
\end{pmatrix}
=\begin{pmatrix}
1 & \frac{1}{2} \\ 0 & \frac{1}{2}
\end{pmatrix}^t

\begin{pmatrix}
2 & 0 \\ 0 & 6
\end{pmatrix}

\begin{pmatrix}
1 & \frac{1}{2} \\ 0 & \frac{1}{2}
\end{pmatrix},

\end{cases}$$
$$(\z_2^{\times})^2=\{u^2|u\in\z_2^{\times}\}=1+8\z_2,$$

$$\left(\z_2^{\times}+\frac{u}{2^r}\right)^2=
\begin{cases}
\frac{1}{2^{r-1}}\z_2^{\times}+\left(\frac{u}{2^r}\right)^2 & \text{when } r>1\\
\frac{1}{2^{r-2}}\z_2^{\times}+\left(\frac{u}{2^r}\right)^2 & \text{when } r=1,
\end{cases}$$
and
$$\left(\z_2+\frac{u}{2^r}\right)^2=
\begin{cases}
\frac{1}{2^{r-1}}\z_2+\left(\frac{u}{2^r}\right)^2 & \text{when } r>1\\
\frac{1}{2^{r-2}}\z_2+\left(\frac{u}{2^r}\right)^2 & \text{when } r=1.
\end{cases}$$
\end{proof}

%%%%%%%%%%%%%%%%%%%%%%%%%%%%%%%%%%%%%%%%%%%55

\section{Globally representing integer which is locally represented}

\begin{lem}\label{main lem}
For an $m$-gonal form $\left<a_1,\cdots,a_n\right>_m$ of rank $n\ge5$, let an integer $A(m-2)+B$ with $0\le B \le m-3$ be locally represented by $\left<a_1,\cdots,a_n\right>_m$.
Then for some 
$$\begin{cases}
k \in [0,K(\mathbf a)]  \\
P=\prod \limits_{p \in T(\mathbf a)\cup\{2\}}p^{s(p)}
\end{cases}$$
with $0\le s(p) \le \frac{1}{2}\ord_p(4a_1)$,
the equation 
\begin{equation} \label{lem eq}
    Q_{a_1;\mathbf a}(P\mathbf x-(B+k(m-2))\mathbf r)=(2A+B+k(m-4))a_1-(B+k(m-2))^2\cdot\left(1-\sum_{i=2}^na_ir_i\right)
\end{equation}
is locally primitively solvable
where $K(\mathbf a)=K(a_1,\cdots,a_n)$ is a constant which is dependent only on $a_1,\cdots,a_n$ and $T(\mathbf a)=T(a_1,\cdots,a_n)$ is a finite set of all odd primes $p$ for which there are at most four units of $\z_p$ in $\{a_1,\cdots,a_n\}$ by admitting a recursion.
\end{lem}
\begin{proof}
By using Propositions \ref{odd,5}, \ref{primitive}, \ref{odd,nd12}, \ref{min ord}, and \ref{2,nd} and Chinese Remainder Theorem, one may prove that this lemma holds for $$K(\mathbf a):=\prod \limits_{p \in {T(\mathbf a)\cup\{2\}}}4\cdot p^{1+\ord_p(a_1)+2\ord_p(d_{a_1;\mathbf a})}-1.$$
\end{proof}

\vskip 0.8em
Now we are ready to prove Theorem \ref{main thm}.
We finish this article by proving Theorem \ref{main thm}.

\vskip 0.8em

\begin{proof}[proof of Theorem \ref{main thm}]
By Theorem 4.9 (2) in \cite{CO}, we may take a constant $$C(P,k;\mathbf a)=C(P,k;a_1,\cdots,a_n)$$ for which 
if 
$$\begin{cases}
\text{$Q_{a_1;\mathbf a}(P\mathbf x-(B+k(m-2))\mathbf r)=c$ is primitively locally solvable} & \text{ and } \\
\text{$c > C(P,k;\mathbf a)$ is sufficiently large,} &
\end{cases}$$
then $ Q_{a_1;\mathbf a}(P\mathbf x-(B+k(m-2))\mathbf r)=c$ is indeed (globally) solvable over $\z$.
Which induces that if \eqref{lem eq} is primitively locally solvable and
$$(2A+B+k(m-4))a_1-(B+k(m-2))^2\cdot\left(1-\sum_{i=2}^na_ir_i\right)>C(P,k;\mathbf a),$$ then $A(m-2)+B$ is represented by the $m$-gonal form $\left<a_1,\cdots,a_n\right>_m$.
Define the constant $C(a_1,\cdots,a_n)$ as the maximum of $C(P,k;\mathbf a)$
where $P$ runs through $\{\prod_{p \in T(\mathbf{a})\cup\{2\}}p^{s(p)}|0\le s(p) \le \frac{1}{2}\ord_p(4a_1) \}$ and $k$ runs through $\{0,1,\cdots, K(\mathbf a)\}$.
Finally we may yield that for an integer $A(m-2)+B$ with $0\le B \le m-3$ which is locally represented by $\left<a_1,\cdots,a_n\right>_m$ if $$A>\frac{1}{2a_1}\left(C(a_1,\cdots,a_n)+(B+K(\mathbf a)(m-2))^2\cdot\left(1-\sum_{i=2}^na_ir_i\right)\right)-\frac{1}{2}(B+k(m-4)),$$
then 
$A(m-2)+B$ is indeed represented by $\left<a_1,\cdots,a_n\right>_m$ by Lemma \ref{main lem}.
And then, the theorem would holds for  $$N(a_1,\cdots,a_n)=\frac{1}{2a_1}\left(C(\mathbf a)+(1+K(\mathbf a))^2\right).$$
\end{proof}

\end{document}